\documentclass[12pt,reqno]{amsart}
\usepackage{amsmath, mathrsfs, dsfont, hyperref, amssymb}
\usepackage{amsthm}
\usepackage[top=1in, bottom=1in, left=1in, right=1in]{geometry}
\usepackage{xcolor}
\usepackage{pdfpages}
\usepackage{graphicx}
\usepackage{caption}
\usepackage{adjustbox}
\usepackage{cite}
\usepackage{lipsum}

\theoremstyle{plain}
\newtheorem{theorem}{Theorem}[section]
\newtheorem{conjecture}[theorem]{Conjecture}

\newtheorem{question}[]{Question}

\begin{document}

\title[Log-concavity And The Multiplicative Properties of Restricted Partition Functions]{Log-concavity And The Multiplicative Properties of Restricted Partition Functions}


\author[Arindam Roy]{Arindam Roy}
\address{University of North Carolina at Charlotte, Department of Mathematics and Statistics, Fretwell 376 9201 University City Blvd.
Charlotte, NC 28223}
\email{aroy15@charlotte.edu}


\begin{abstract}
    The partition function $p(n)$ and many of its related restricted partition functions have recently been shown independently to satisfy log-concavity: $p(n)^2 \geq p(n-1)p(n+1)$ for $n\geq 26$, and satisfy the inequality: $p(n)p(m) \geq p(n+m)$ for $n\geq m\geq 2$ with only finitely many instances of equality or failure. This paper proves that this is no coincidence, that any log-concave sequence $\{x_n\}$ satisfying a particular initial condition likewise satisfies the inequality $x_nx_m \geq x_{n+m}$. This paper further determines that these conditions are sufficient but not necessary and considers various examples to illuminate the situation. 
\end{abstract}

\maketitle
\section{Introduction}
The \textit{partition function} $p(n)$ enumerates the number of partitions of a positive integer $n$ where the partitions are positive integer sequences $\lambda=(\lambda_1,\lambda_2,...)$ with $0<\lambda_1\leq\lambda_2\leq~\dots$ and $\sum_{j\geq1}{\lambda_j}=n$. For example, $p(4)=5$ since the only ways to partition $4$ are $4$, $3+1$, $2+2$, $2+1+1$, and $1+1+1+1$. Recent interest regarding the various properties of $p(n)$ led to a conjecture by Chen \cite{[Chen]}, and subsequent proof by DeSalvo \& Pak \cite{[Salvo]}, that the sequence $\{p(n)\}_{n\in\mathbb{N}}$ is log-concave for $n>25$.  A sequence of integers $\{x_n\}_{n\in\mathbb{N}}$ is \textit{log-concave} if $x_n^2~\geq~x_{n-1}x_{n+1}$ and is \textit{log-convex} if $x_n^2 \leq x_{n-1}x_{n+1}$ for all $n\geq n_0$ for some $n_0\in \mathbb{N}$. 

Log-concavity had been previously established in 1978 by Nicholas \cite{[Nicolas]} for the partition function, but there was a revitalization of interest following the result of DeSalvo \& Pak \cite{[Salvo]}. New results found in the work of Bessenrodt \& Ono \cite{BessOno} in 2014 proved that the sequence$\{p(n)\}_{n\in\mathbb{N}}$ satisfies the inequality
\begin{align}\label{bres-ono}
    p(n)p(m) \geq p(n+m)\quad\text{for all}\quad (m,n)\quad\text{with}\quad n\geq m\geq 2
\end{align}
  where instances of equality occur at $(2,6), (2,7), (3,4)$, and instances of failure occur at $(2,2)$, $(2,3)$, $(2,4)$, $(2,5)$, $(3,3)$, $(3,5)$ and nowhere else. 

For convenience, this paper uses the notation: a sequence $\{x_n\}$ is \textit{multiplicatively abundant} (or simply \textit{abundant}) if $x_nx_m \geq x_{n+m}$, and is \textit{multiplicatively deficient} (or simply \textit{deficient}) if $x_nx_m \leq x_{n+m}$ for all $n\geq m\geq n_0$ but finitely many $n\geq m\geq n_0$. In these terms, Bessenrodt \& Ono \cite{BessOno} proved that $\{p(n)\}$ is abundant for all $n\geq m\geq 2$. Proofs that the partition function is log-concave or abundant initially relied on an analytic result of Rademacher type due to Lehmer \cite{[Lehmer]} until 2017 when A. Alanazi, S.M. Gagola III, \& A.O. Munagi \cite{[Alanzi]} published a combinatorial proof that $\{p(n)\}$ is abundant.

Many sequences (especially sequences arising from restricted partitions) are sporadic for small indices but become and stay log-concave for sufficiently large indices. The interaction between log-behavior and multiplicative behavior in the setting of general sequences was studied in the early 2000's by Asai, Kubo, \& Kuo \cite{[Asai]} 
who collectively established inequalities for Bell numbers, Motzkin numbers, Fine numbers, Fresnel numbers, Ap\`ery numbers, large Schr\"oder numbers, central Delannoy numbers, and Baxter permutations. Recently, Gajdzica, Miska, \& Ulas \cite{Gajdzica_Miska_Ulas} studied similar interactions in a more general setting. 

The main goal of this paper is to prove that \textit{any} sequence that is eventually log-concave is also eventually abundant, provided the sequence also satisfies a sufficient initial conditions. 
\begin{theorem}\label{main_theorem1}
Let $\{x_n\}$ be a sequence of positive real numbers and let $N$ be a non-negative integer such that for all $n > N $, the sequence satisfies $x_n^2 \geq x_{n-1}x_{n+1}$. Then for all $n >N $ and $m\geq 0$, the sequence $\{x_n\}$ satisfies:
\begin{align}
(x_{N})^{\frac{m}{n-N+m}}(x_{n+m})^{\frac{n-N}{n-N+m}}\leq x_n \leq x_{N}\left(\frac{x_{N+1}}{x_{N}}\right)^{n-N}. \label{boundsineq}
\end{align}
Moreover, if there exists an integer $k\geq 0$ such that
\begin{align}\label{Condition 1}
    (x_{N+k})^{N+k+1}\geq (x_{N+k+1})^{N+k}
\end{align} 
then $x_{n+m} \leq x_nx_m$ for all $n,m \geq N+k $. In addition, $\{x_n\}$ is abundant for all $m,n\geq d$ with at most finite number of failures $(m,n)\in [d,M+1]\times [d,N+k]\bigcup [d,N+k]\times [d,M+1]$, if $d\leq N+k$ is the smallest positive integer such that \begin{align}\label{condition extra} \left(\frac{x_{M+1}}{x_M}\right)^m\leq x_m \end{align} for all $d\leq m\leq N+k$ and for some $M\geq N+k$.
\end{theorem}

\noindent {\bf{Remarks:}}
1) If a log-concave sequence does not satisfy condition \eqref{Condition 1}, i.e., $(x_{N+k})^{N+k+1}\geq (x_{N+k+1})^{N+k}$ for all $k\geq 0$, then it is not true that $x_{n+m} \leq x_nx_m$. For example,
define the sequence $\{a_n\}$ by $a_n=2^{n}-1/2$ for $n = 0,1,2,\ldots$
The sequence is \textit{strictly} log-concave: $$
\left(2^{n}-\frac{1}{2}\right)^2=2^{2n}-2^{n}+\frac{1}{4} > 2^{2n} - 2^{n} - 2^{n-2} +\frac{1}{4}=\left(2^{n+1}-\frac{1}{2}\right)\left(2^{n-1}-\frac{1}{2}\right).$$
However, the sequence $\{(2^{n}-1/2)^{\frac{1}{n}}\}$ for $n = 1,2, 3,  \ldots$, is a strictly increasing sequence. 
Hence, $a(n)$ does not satisfy Condition \eqref{Condition 1}. Indeed,
$
(2^{n+k}-1/2)^{n+k+1} < (2^{n+k+1}-1/2)^{n+k} 
$ 
for all $n, k\geq 0$. In this particular case, $\{a_n\}$ is also not abundant (in fact, $\{a_n\}$ is deficient):
\[
(2^{n}-1/2)(2^{m}-1/2) = 2^{n+m}-2^{n-1}-2^{m-1}+\frac{1}{4} < 2^{n+m}-1/2
\]
for $n\geq m \geq 0$. 
\\

\noindent 2) Suppose $N=0$, i.e.,  $\{x_n\}$ is log-concave for all $n\geq 1$,  then 
\[
x_0^{\frac{m}{n+m}}(x_{n+m})^{\frac{n}{n+m}} \leq x_n \leq x_{0}\left(\frac{x_{1}}{x_{0}}\right)^{n}
\]
for all non-negative integers $m$. This give an upper bound for a log-concave sequence $\{x_n\}$. For example, a sequence with super exponential growth can not be log-concave. Moreover, $x_{n+m}\leq x_nx_m$ for all $n\geq m\geq k$ if $(x_{k})^{k+1}\geq (x_{k+1})^{k}$ for any fixed non-negative integer $k$. Specifically, if $k=0$, meaning that $x_0\geq 1$, then for all $n\geq m\geq 0$, $x_{n+m}\leq x_nx_m$. In 2001, Asai, Kubo, \& Kuo \cite[Theorem 2]{[Asai]}\footnote{While $x_0=1$ is mentioned in the theorem's statement, the proof works for $x_0\leq 1$. To obtain the proof for the log-concave sequence, one must take the reciprocal of the log-convex sequence in \cite{[Asai]}.} proved the case $N=0$ and $k=0$ of Theorem \ref{main_theorem1}. 
Gajdzica, Miska, \& Ulas \cite[Theorem 4.4]{Gajdzica_Miska_Ulas} proved this case again in $2023$ using a similar technique used in \cite{[Asai]}. It can be seen that Theorem \ref{main_theorem1} not only generalized the results of \cite[Theorem 2]{[Asai]} or \cite[Theorem 4.4]{Gajdzica_Miska_Ulas}, it is also more effective. For instance, Bessenrodt-Ono result \eqref{bres-ono} can be obtained from Theorem \ref{main_theorem1} but it can not be derived from the $N=k=0$ case of Theorem \ref{main_theorem1}. Instead one may obtain a weaker result such as $\{p(n+26)\}_{n\geq 0}$ is abundant (see \cite[Example 4.5]{Gajdzica_Miska_Ulas}). A weaker version of Theorem \ref{main_theorem1} for $k=0$ and $M=N$ was recently derived by Heim and Neuhauser \cite{heimmarkus}. As we shall demonstrate in section 1.1, the result given in \cite{heimmarkus} does not imply that the log-concave sequences in example 5 are abundant for smaller indices.

\subsection{Some Examples of Theorem \ref{main_theorem1}}
1) 
In light of Theorem \ref{main_theorem1} and because log-concavity has been established for $n>25$ in \cite{[Salvo], [Nicolas]}, to show the partition function is abundant it suffices to show that Condition \eqref{Condition 1} holds for $p(25)=1958$ and $p(26)=2436$ with $N=25$ and $k=0$: it is easy to verify that in fact, $1958^{\frac{1}{25}} > 2436^{\frac{1}{26}}$. 
It can be noted that this result only proves the partition function is abundant for $n\geq m\geq 25$ but not the result of Bessenrodt-Ono given in \eqref{bres-ono}. To show \eqref{bres-ono} holds one needs to check condition \eqref{condition extra}, in particular, for $M=25$: $p(26)/p(25)\leq (p(m))^{1/m}$ for all $2\leq m\leq 25$. This can be easily verified. Next, by checking in the region $ [2,26]\times [2,26]$, we see failure occur at $(2,2)$, $(2,3)$, $(2,4)$, $(2,5)$, $(3,3)$, $(3,5)$, $(3,2)$, $(4,2)$, $(5,2)$, $(5,3)$ and equality occur at $(2,6), (2,7), (3,4), (6,2), (7,2), (4,3)$. This gives the result of Bessenrodt-Ono \eqref{bres-ono}.
\\

\noindent 2) 
The \textit{overpartition function}, $\overline{p}(n)$, counts the number of partitions of $n$ in which the last occurrence of a number can be distinguished.  
For example, $\overline{p}(4)=14$: 
\begin{align*}
&\ 1+1+1+1,\quad 1+1+1+\overline{1},\quad 1+1+2,\quad 1+1+\overline{2},\quad 1+\overline{1}+2,\\
&1+\overline{1}+\overline{2},\quad 2+2,\quad 2+\overline{2},\quad 1+3,\quad 1+\overline{3},\quad \overline{1}+3,\quad \overline{1}+\overline{3},\quad 4,\quad \overline{4}.
\end{align*}
\[
    \overline{p}(n)=1, 2, 4, 8, 14, 24, 40, 64, 100, 154, 232, 344, 504, 728, 1040 ,\ldots \text{OEIS: A015128 \cite{[A015128]}}
    \]
Combining the result of Engel \cite{[Engel]} and the fact $\overline{p}^2(1)\geq \overline{p}(0)\overline{p}(2)$ we see that $\{\overline{p}(n)\}$ is log-concave for $n>0$. Note that $\overline{p}$ satisfies both Conditions \eqref{Condition 1} and \eqref{condition extra} for $N=k=0$. Then Theorem \ref{main_theorem1}  give us that the overpartition function satisfies the inequality 
\begin{align}
\overline{p}(m)\overline{p}(n)\geq \overline{p}(m+n)\quad\text{for all}\quad  m,n\geq 0 \label{liuzhang}
\end{align}
     where instances of equality occur at $(1,1)$, $(1,2)$, $(2,1)$, $(m,0)$ and $(0,n)$.
Note that \eqref{liuzhang} was proven using the asymptotic of the overpartition function by Liu \& Zhang \cite{liu} in 2021.
\\

\noindent 3) The \textit{distinct Partition Function} $P_d(n)$ counts the number of ways to partition $n$ into distinct parts, where for all $i$ and $j$, $\lambda_i\neq\lambda_j$.
\begin{align*}
P_d(n)&=1, 1, 1, 2, 2, 3, 4, 5, 6, 8, 10, 12, 15, 18, 22, 27, 32, 38, 46, 54, 64, 76, 89, 104, 122, 142, 165,\\ &\quad 192, 222, 256, 296, 340, 390, 448, 512, 585, 668,\ldots \text{OEIS: A000009 \cite{[000009]}}
\end{align*}
Curiously, Euler showed that $P_d(n)$ is equivalent to the number of partitions of $n$ using only odd integers \cite{[Euler],[Kim],[Merca]}. For example, $P_d(4)=2$ since there are only two ways to partition $4$ into distinct parts: $4=4=3+1$ (and only two ways to partition $4$ into odd parts: $4=3+1=1+1+1+1$). In 2023, Dong \& Ji \cite{[Dong]} proved that the sequence $\{P_d(n)\}$ is log-concave for all $n>32$. Note that $P_d(n)$ satisfies Condition \eqref{Condition 1} for $N=32$ and $k=0$: $390^{\frac{1}{32}}>448^{\frac{1}{33}}$ and condition \eqref{condition extra} for $M=N=32$: $448/390\leq p_d(m)^{1/m}$ for all $3\leq m\leq 32$. Then from Theorem \ref{main_theorem1} and numerical computation gives the distinct partition function satisfies the inequality 
\begin{align}
    P_d(m)P_d(n) \geq P_d(m+n)\quad\text{for all}\quad n\geq n>2 \label{distinct_partition}
\end{align}
where instances of equality and failure occur as given in Table \ref{P_d(n)_table}.
\begin{table}[ht]
\caption{Equality and Failure for $P_d(m)P_d(n)\geq P_d(m+n)$}
\begin{tabular}{c|c}\label{P_d(n)_table}
Equality at $(m,n)$&Failure at $(m,n)$\\
\hline
$(1,n)$; $n=1,3$& $(1,n), (2,n)$ ; $n\geq 2$  \\
$(3,n)$; $n=3,5,6,7,8$& $(1,2)$, $(3,4)$, $ (5,5)$\\
$(4,n)$; $n=15,16,17$& $(4,n)$; $4\leq n \leq 14,$\\
$(5,n)$; $n=6,7,8$ & 
\end{tabular}
\end{table}
The inequality \eqref{distinct_partition} together with the table \ref{P_d(n)_table} was proven by Beckwith \& Bessenrodt in \cite{[Beckwith]} using a Rademacher type asymptotic formula of $k$-regular partition.
\\

\noindent 4) This is an example of $k>0$ case. 
Define the sequences $a_0=\frac{1}{2}$ and $a_n=2^{\sqrt{n}}$ for $n\geq 1$. Keep in mind that $a_1^2=2^2>\frac{1}{2}\times 2^{\sqrt{2}}=a_0a_2$, and that $a_{n+1}^2=2^{2\sqrt{n+1}}>2^{\sqrt{n}}\times 2^{\sqrt{n+2}}=a_na_{n+2}$ for all $n\geq 1$ and that $a_{n+m}=2^{\sqrt{n+m}}<2^{\sqrt{n}}\times 2^{\sqrt{m}}=a_na_m$ for all $n,m\geq 1$.  In other words, $\{a_n\}$ fulfills the inequality $a_{n+m}<a_na_m$ for all $n,m\geq 1$, but is a log-concave sequence for any $n$ with $a_0<1$. Since $a_0<1$, the result in \cite[Theorem 4.4]{Gajdzica_Miska_Ulas} is not applicable here. 
It is hopefully clear that $(a_1)^2=2^2>2^{\sqrt{2}}=(a_2)^1$ which satisfy Condition \eqref{Condition 1} in Theorem \ref{main_theorem1}.
\\

\noindent 5) Lastly, we give an example when $M>N+k$. Consider the sequence $a_n=\left(\frac{31}{25}\right)^n+\frac{1}{10^6}$ when $0\leq n\leq 24$ and $a_n=p(n)$ when $n\geq 25$. Here $p(n)$ is the partition function. It can be checked that $a_n^2\geq a_{n-1}a_{n-2}$ for $n\geq 25$. Also, it can be verified that $\displaystyle(a_n)^{1/n}<\frac{a_{26}}{a_{25}}<\frac{a_{25}}{a_{24}}$ for all $1\leq n\leq 24$ but $\displaystyle(a_n)^{1/n}>\frac{a_{27}}{a_{26}}$ for all $1\leq n\leq 25$. Thus, using Theorem \ref{main_theorem1}, we can deduce that, for every $m,n\geq 1$, $a_{n+m}\leq a_na_m$, with the exception of finitely many failures in $[1,25]\times [1,25]$. Note that, the result in \cite{heimmarkus} yields a weaker result, which only gives $a_{n+m}\leq a_na_m$ for all $m,n\geq 24$.
\\
\\

Consider what happens to the conditions of Theorem \ref{main_theorem1} if the sequence $x_n$ is replaced with its reciprocal sequence. The analogue to Condition \eqref{Condition 1} for log-convexity is similarly sufficient to ensure that a sequence is deficient. Consider the following equivalent formulation of Theorem \ref{main_theorem1} in terms of log-convexity:
\begin{theorem}\label{main_theorem2}
Let $\{y_n\}$ be a sequence of positive real numbers and let $N$ be a positive integer such that for all $n > N $, the sequence satisfies $y_n^2 \leq y_{n-1}y_{n+1}$. Then for all $n > N $ and $m\geq0$, the sequence $\{y_n\}$ satisfies:
\begin{align*}
\left(y_{N}\right)^{\frac{m}{n-N+m}}\left(y_{n+m}\right)^{\frac{n-N}{n-N+m}} \geq y_n \geq y_{N}\left(\frac{y_{N+1}}{y_{N}}\right)^{n-N}.
\end{align*}
Moreover, if there exists an integer $k\geq 0$ such that
\begin{align}
    (y_{N+k})^{N+k+1}\leq (y_{N+k+1})^{N+k},\label{Condition 2}
\end{align}
then $y_{n+m} \geq y_ny_m$ for all $n,m \geq N+k $. In addition, $\{y_n\}$ is deficient for all $m,n\geq d$ with at most finite number of failures $(m,n)\in [d,M+1]\times [d,N+k]\bigcup [d,N+k]\times [d,M+1]$, if $d\leq N+k$ is the smallest positive integer such that \begin{align}\label{condition extra1} \left(\frac{y_{M+1}}{y_M}\right)^m\geq y_m \end{align} for all $d\leq m\leq N+k$ and for some $M\geq N+k$.
\end{theorem}
\subsection{Sufficient but not Necessary}
\indent Theorem \ref{main_theorem1} establishes the condition of log-concavity together with Condition \eqref{Condition 1} as sufficient conditions that a sequence is abundant, but are they necessary conditions? In fact, there are sequences that are not log-concave (nor log-convex) yet are abundant.
\\
{\bf Example 1.}\label{Example 3}
Define the sequence $a_n$ by
\[
a_n=\begin{cases}
		2, & \text{if}\ n\equiv0 \mod{4}\\
            3, & \text{if}\ n\equiv1,3\mod{4}\\
            4, & \text{if}\ n\equiv2 \mod{4}.
\end{cases}
\]
\\
The inequality $a_n^2 \geq a_{n-1}a_{n+1}$ fails whenever $n\equiv0\mod{4}$, and so $\{a_n\}$ is never log-concave. Similarly, the inequality $a_n^2 \leq a_{n-1}a_{n+1}$ fails whenever $n\equiv2\mod{4}$, and so $\{a_n\}$ is never log-convex. However, the sequence $\{a_n\}$ satisfies the inequality $a_na_m \geq a_{n+m}$ for all $n$ and $m$. For $n>0$, note that $(a_{n+1})^\frac{1}{n+1} \leq (a_{n+2})^\frac{1}{n+2}$ whenever $n\equiv0,3 \mod{4}$ and $(a_{n+1})^\frac{1}{n+1} \geq (a_{n+2})^\frac{1}{n+2}$ whenever $n\equiv1,2 \mod{4}$.
\\
\indent Analogously, Theorem \ref{main_theorem2} establishes the condition of log-convexity together with Condition \eqref{Condition 2} as sufficient conditions that a sequence is deficient, but are they necessary conditions? In fact, there are sequences that are neither log-concave nor log-convex yet are deficient. \\
\\
{\bf Example 2.}\label{Example 4} The Fibonacci sequence defined by $F_0=0$, $F_1=1$, and $F_n=F_{n-1}+F_{n-2}$, is log-concave for odd indices and log-convex for even indices (known as \textit{log-Fibonacci} \cite{[Asai]}). This is demonstrated in the literature by Cassini's identity \cite{[Spicey]}:
\[
F_{n+1}F_{n-1}-F_n^2=(-1)^n
\]  for all $n$. 
The Fibonacci sequence is also known \cite{[Honsberger]} to satisfy the identity
\[
F_{n+m}=F_{m+1}F_n+F_mF_{n-1}
\]
for all $n$ and $m$. Note that $F_{m+1}F_n+F_mF_{n-1} > F_nF_m$. Then, besides the single instance of equality at $(n,m)=(1,1)$, the Fibonacci sequence satisfies the inequality $F_{n+m} > F_nF_m$ for all $n$ and $m$. Note that for all $n>1$, the sequence $\{(F_{n})^{\frac{1}{n}}\}$ is a strictly increasing sequence that approaches the golden ratio in the limit, hence $(F_{n+1})^{\frac{1}{n+1}} < (F_{n+2})^{\frac{1}{n+2}}$.
\subsection{Other Log-Concave Sequences} There are many examples of log-concave sequences that, in light of Theorem \ref{main_theorem1}, are now known to be abundant. Log-concavity has been established for Stirling numbers of the first and second kind \cite{[SIBUY]}, the Hyperfibonacci and Hyperlucas numbers \cite{[Zheng]}, the HyperPell and HyperPell-Lucas numbers \cite{[Ahmiaab_Belbachirb_Belkhirb]}, and the figurate and central figurate numbers \cite{[Gedefa_1],[Gedefa_2]} to name a few. A more general class of sequences is proven log-concave by O'Sullivan \cite{[Sullivan2]} who considered sequences $\{a(n)\}$ of the form $a(n)=n^r\exp(tn^s)$ for various values of $r,s$ and $t$.

It is natural to consider  \textit{restricted partition functions}, frequently denoted $p_A(n)$, where $A$ is some parameter limiting which values $\lambda$ can take. For instance, the integer $4$ could be partitioned among the restricted set of prime numbers, or perfect squares, or powers of $2$, or odd numbers, etc:
\begin{alignat*}
4  & 4= 2+2  &\ \ \qquad  4 &  = 2^2            &\ \ \qquad  4  & = 2^2  &\ \ \qquad  4  & = 3+1\\
   &        &          & = 1^2+1^2+1^2+1^2 &           & = 2^1+2^1  &                       &=1+1+1+1\\
  &        &          &                      &           & = 2^0+2^0+2^0+2^0  &         
\end{alignat*}

Proofs that particular restricted partition functions are log-concave or abundant or both have typically relied on knowing their asymptotic. The $k$-regular partition function for $2\leq k\leq 6$ was proven abundant in 2016 by Beckwith \& Bessenrodt \cite{[Beckwith]}, and log-concavity was proven in 2019 by Craig \& Pun \cite{[Craig]} for large $n$. Proofs of log-concavity for the $k$-regular partition with $k=2$ and $3$ are given by Dong \& Ji \cite{[Dong]}. The $k$-colored partition function was proven abundant in 2018 by Chern, Fu, \& Tang \cite{[Chern]}. Log-concavity of the $k$-colored partition was proved by Bringmann, Kane, Rolen, \& Tripp \cite{BKRT} which was conjectured in \cite{[Chern]}. Dawsey \& Masri \cite{[Dawsey]} proved the Andrews smallest parts function is log-concave and abundant in 2017. In 2022, the plane partition function was proven abundant by Heim, Neuhauser, \& Tröger \cite{[Heim]} who also conjectured log-concavity; this was later proven that same year by Ono, Pujahari, \& Rolen \cite{[Onoetal]}.

There are restricted partition functions for which log-concavity has been established, but which have not been proven to be abundant. Proofs of log-concavity for such functions are given 
by Iskander, Jain, \& Talvola \cite{[Iskander]} for the \textit{fractional partition function}, by Gajdzica \cite{[Gajdzica]} for the \textit{multiset partition function}, and  independently by Benfield, Roy, \& Paul \cite{benfield} and by O'Sullivan \cite{[Sullivan]} for the \textit{$k$-th power partition function}. 

This paper proves that these restricted partition functions are abundant. Where possible for these restricted partition functions, in the spirit of the original publication by Bessenrodt and Ono \cite{BessOno}, this paper also calculates the finitely many values where the strict abundant inequality fails or is equal.
\section{Proof of Theorem \ref{main_theorem1}}

Let us  assume $x_k^2 \geq x_{k-1}x_{k+1}$ where $k=N+1,N+2,\dots,n$. Now consider the two inequalities $x_n^2 \geq x_{n-1}x_{n+1}$ and $x_{n-1}^2 \geq x_nx_{n-2}$. After some rearrangement, it follows that 
\begin{align}\label{neq1}
x_{n+1} \leq \frac{(x_n)^2}{x_{n-1}} \leq \frac{1}{x_{n-1}}\left( \frac{(x_{n-1})^2}{x_{n-2}}\right)^2 = \frac{(x_{n-1})^3}{(x_{n-2})^2}.
\end{align}
Since $x_{n-2}^2 \geq x_{n-1}x_{n-3}$,  substituting this into \eqref{neq1}, yielding 
\begin{align}
\nonumber \frac{(x_{n-1})^3}{(x_{n-2})^2} \leq \frac{1}{(x_{n-2})^2}\left( \frac{(x_{n-2})^2}{x_{n-3}} \right)^3 = \frac{(x_{n-2})^4}{(x_{n-3})^3}.
\end{align}
Iterating this process 
yields the following chain of inequalities:
\begin{align}\label{observation_1}
x_{n+1} \leq \frac{(x_{n-1})^3}{(x_{n-2})^2} \leq \frac{(x_{n-2})^4}{(x_{n-3})^3} \leq 
\frac{(x_{n-3})^5}{(x_{n-4})^4} \leq \ldots \leq \frac{(x_{N+2})^{n-N}}{(x_{N+1})^{n-N-1}} \leq x_{N}\left(\frac{x_{N+1}}{x_{N}}\right)^{n+1-N}.
\end{align}
Hence, $x_n$ is bounded from above: $x_n \leq x_{N}\left(\frac{x_{N+1}}{x_{N}}\right)^{n-N}$ for all $n>N$ . This gives the right-hand side of \eqref{boundsineq}. 

To prove the left-hand side of \eqref{boundsineq}, note that the inequalities $(x_{N+1})^2 \geq x_{N}x_{N+2}$ and $(x_{N+2})^2 \geq x_{N+1}x_{N+3}$ follow from log-concavity. After some substitution, these inequalities can be written as
\begin{align}\label{blbl}
    (x_{N+2})^2 \geq x_{N+1}x_{N+3}  \geq \left((x_{N})^\frac{1}{2}(x_{N+2})^\frac{1}{2}\right)x_{N+3} \Rightarrow x_{N+2} \geq (x_{N})^\frac{1}{3}(x_{N+3})^\frac{2}{3}.
\end{align}
Again by log-concavity, $(x_{N+3})^2 \geq x_{N+2}x_{N+4}$ which after substitution into \eqref{blbl} can be written as
\begin{align}
    \nonumber(x_{N+3})^2 &\geq x_{N+2}x_{N+4}  \geq \left((x_{N})^\frac{1}{3}(x_{N+3})^\frac{2}{3}\right)x_{N+4} \Rightarrow x_{N+3} \geq (x_{N})^\frac{1}{4}(x_{N+4})^\frac{3}{4}.
\end{align}
This process may continue recursively, and as $n$ gets large, 
\begin{align}
   x_n \geq \left(x_{N}\right)^\frac{1}{n-N+1}\left(x_{n+1}\right)^\frac{n-N}{n-N+1}. \label{eq4}
\end{align}
 Increasing the index from $n$ to $n+1$ yields
\begin{align*}
    x_{n+1} \geq \left(x_{N}\right)^\frac{1}{n-N+2}\left(x_{n+2}\right)^\frac{n-N+1}{n-N+2}
\end{align*}
which can be substituted into \eqref{eq4} and rearranged:
\begin{align*}
    \nonumber x_n \geq \left(x_{N}\right)^\frac{1}{n-N+1}\left(x_{n+1}\right)^\frac{n-N}{n-N+1} &\geq \left(x_{N}\right)^\frac{1}{n-N+1}\left(\left(x_{N}\right)^\frac{1}{n-N+2}\left(x_{n+2}\right)^\frac{n-N+1}{n-N+2}\right)^{\frac{n-N}{n-N+1}}\\
   &= (x_{N})^{\frac{1}{n-N+1}+\frac{n-N}{(n-N+1)(n-N+2)}}(x_{n+2})^\frac{n-N}{n-N+2}.
\end{align*}

Iterating this process up to some index $n+m$ for some $m\geq 0$ yields:
\begin{align}
    x_n \geq (x_{N})^{\frac{1}{n-N+1}+\frac{n-N}{(n-N+1)(n-N+2)}+\dots+\frac{n-N}{(n-N+m-1)(n-N+m)}}(x_{n+m})^\frac{n-N}{n-N+m}.\label{exponent}
\end{align}

Note that inequality \eqref{exponent} holds for all $n>N$ and all $m\geq 0$, and that the exponent of $x_{N}$ in \eqref{exponent} is a telescoping series: 
\begin{align*}
    &\frac{1}{n-N+1}+\frac{n-N}{(n-N+1)(n-N+2)}+\dots+\frac{n-N}{(n-N+m-1)(n-N+m)}\\
    &=(n-N)\left[\left(\frac{1}{n-N}-\frac{1}{n-N+1}\right)+\left(\frac{1}{n-N+1}-\frac{1}{n-N+2}\right)+\ldots\right.\\
    &\hspace{3in}\left.+\left(\frac{1}{n-N+m-1}-\frac{1}{n-N+m}\right)\right]\\
    &=(n-N)\left(\frac{1}{n-N}-\frac{1}{n-N+m}\right)\\
    &=\frac{m}{n-N+m}.
\end{align*}
Thus, inequality \eqref{exponent} is reduced to:
\begin{align}\label{observation_2}
    x_{n} \geq (x_{N})^{\frac{m}{n-N+m}}(x_{n+m})^{\frac{n-N}{n-N+m}}
\end{align}
 for all $n> N$. This gives the left-hand side of \eqref{boundsineq}.

To prove the abundantness, it is necessary to establish relations for all entries in the sequence. 
To achieve this, construct a sequence $\{a_n\}$ 
recursively by $a_n=x_{n}$ for all $n \geq N+k$ and for $n< N+k$, 
\begin{align}
    &a_{N+k-1} := \frac{(a_{N+k})^2}{a_{N+k+1}} = \frac{(x_{N+k})^2}{x_{N+k+1}}\\\nonumber
    &a_{N+k-2} := \frac{(a_{N+k-1})^2}{a_{N+k}} = \frac{\left(\frac{(x_{N+k})^2}{x_{N+k+1}}\right)^2}{x_{N+k}} = \frac{(x_{N+k})^3}{(x_{N+k+1})^2}\\\nonumber
    &\vdots\\
    &a_0 = \frac{(x_{N+k})^{N+k+1}}{(x_{N+k+1})^{N+k}}\label{a0def}.
\end{align}
Hence, $(a_{n+1})^2\geq a_na_{n+2}$ for all $n\geq 0$ and $a_0\geq 1$ by Condition \eqref{Condition 1}. Now from Asai, Kubo, \& Kuo \cite[Theorem 2]{[Asai]} or Gajdzica, Miska, \& Ulas \cite[Theorem 4.4]{Gajdzica_Miska_Ulas}, one obtains \begin{align}a_na_m\geq a_{n+m} \quad\text{for all}\quad m,n\geq 0.\label{abun} \end{align}

Here, we present an alternative proof of \eqref{abun} compared to the one given in \cite{[Asai],Gajdzica_Miska_Ulas}, as it aligns better with the approach used in the initial part of the proof.
A similar process to that outlined in \eqref{observation_1} yields
\begin{align}
a_{n+1} \leq \frac{(a_{n-1})^3}{(a_{n-2})^2} \leq \frac{(a_{n-2})^4}{(a_{n-3})^3} \leq 
\frac{(a_{n-3})^5}{(a_{n-4})^4} \leq \ldots \leq \frac{(a_2)^n}{(a_1)^{n-1}}=a_1\left(\frac{a_2}{a_1}\right)^{n}.\label{anba1}
\end{align}
Also, by the same process outlined between equations \eqref{blbl} and \eqref{observation_2} one yields
\begin{align}\label{observation_3}
    a_{n} \geq (a_{1})^{\frac{m}{n+m-1}}(a_{m+n})^{\frac{n-1}{n+m-1}}.
\end{align}
Switching the indices $n$ and $m$ in \eqref{observation_3} gives the inequality:
\begin{align}\label{observation_2_swap_m_&_n}
    a_{m} \geq (a_{1})^{\frac{n}{m+n-1}}(a_{m+n})^{\frac{m-1}{m+n-1}}.
\end{align}
Suppose for contradiction that the  sequence $\{a_n\}$ is log-concave but  the sequence $\{a_n\}$ does not satisfies the inequality $a_{n+m} \leq a_na_m$ for all $n,m\geq 0$. Then, for some $m\geq 0$ and $n\geq 0$, the inequality $a_{n+m}>a_na_m$ holds. Hence, for one such pair of $m$ and $n$,  and from \eqref{observation_2} and \eqref{observation_2_swap_m_&_n} it follows that
\begin{align*}
    a_{n+m}&>a_{n}a_{m}\\ &\geq \left((a_{1})^{\frac{m}{n+m-1}}(a_{n+m})^{\frac{n-1}{n+m-1}}\right)\left((a_{1})^{\frac{n}{m+n-1}}(a_{m+n})^{\frac{m-1}{m+n-1}}\right)\\
    &=(a_{1})^{\frac{n+m}{n+m-1}}(a_{n+m})^{\frac{n+m-2}{n+m-1}}.
\end{align*}
Hence
\begin{align}\label{contradiction_left}
    (a_{n+m})^{\frac{1}{n+m-1}} > (a_{1})^{\frac{1}{n+m-1}+1}.
\end{align}
Note by Conditions \eqref{Condition 1} and \eqref{a0def} that $a_0\geq 1$. And by log-concavity, $a_{1}^2 \geq a_{0}a_{2} \geq a_{2}$ implying that $ \frac{a_{2}}{a_{1}}\leq a_{1}$. Then from \eqref{anba1} it follows that
\begin{align}
    \nonumber (a_{n+m})^{\frac{1}{n+m-1}} \leq \left[a_{1}\left(\frac{a_{2}}{a_{1}}\right)^{n+m-1}\right]^{\frac{1}{n+m-1}} &\leq \left(\frac{a_{2}}{a_{1}}\right)\left(a_{1}\right)^{\frac{1}{n+m-1}}\\
    &\leq (a_{1})^{\frac{1}{n+m-1}+1}.\label{contradiction_right}
\end{align}
Hence both \eqref{contradiction_left} and \eqref{contradiction_right} cannot be simultaneously true, this gives a contradiction.
Hence $a_{m+n}\leq a_ma_n$ for all $m,n\geq 0$ and in particular $m,n\geq N+k$.   

Now, since $M\geq N+k\geq N$ then following the same steps of the proof of \eqref{observation_2} we find 
\begin{align*}
    x_{n} \geq (x_{M})^{\frac{m}{n-M+m}}(x_{n+m})^{\frac{n-M}{n-M+m}}
\end{align*}
and in particular 
\begin{align*}
    x_{M+1} \geq (x_{M})^{\frac{m}{m+1}}(x_{M+1+m})^{\frac{1}{m+1}}.
\end{align*}
This gives 
\begin{align*}
    \frac{x_{M+1+m}}{x_{M+1}}\leq \left(\frac{x_{M+1}}{x_{M}}\right)^m.
\end{align*}
Therefore by Condition \eqref{condition extra}
\begin{align*}
    \frac{x_{M+1+m}}{x_{M+1}}\leq x_m
\end{align*}
for all $d\leq m\leq N+k$. Since the sequence $\{\frac{x_{n+1}}{x_n}\}_{n\geq N}$ is decreasing then
\begin{align*}
    \frac{x_{n+m}}{x_{n}}\leq\frac{x_{M+1+m}}{x_{M+1}}\leq x_m
\end{align*}
for all $n\geq M+1$ and $d\leq m\leq N+k$. This completes the proof of Theorem \ref{main_theorem1}.


\section{Various Partition Functions}
This section considers restricted partition functions $P_A(n)$ that are known to be log-concave. It follows from Theorem \ref{main_theorem1} that these sequences are also abundant, provided they also satisfy Condition \eqref{Condition 1}. This section considers a few, well-studied sequences and partition functions and determines for which small indices abundance fails. 
The following subsections have been proven to satisfy log-concavity for sufficiently large $n$, however, the smallest such index $N$ such that each sequence is log-concave for all $n \geq N$ is either unknown or only conjectured. On the other hand, it is known that the asymptotics of these sequences are each of the form $n^re^{t\cdot n^s}$ for real numbers $s, r,$ and $t$ where $0 < s \leq 1$. This implies that for some $n\geq n_0$ the $n$\textsuperscript{th} root of these sequences are each decreasing, hence, there is an index $N$ such that $x_{N+1}^{\frac{1}{N+1}} \geq x_{N+2}^{\frac{1}{N+2}}$.

\subsection{The \texorpdfstring{$k$}\ -th Power Partition Function}\ 
\newline
The \textit{$k$-th power partition function} $p^k(n)$ counts the number of partitions of $n$ among perfect $k$-th powers. For example, $p^k(n)=p(n)$ when $k=1$, and for $k=2$, $p^2(n)$ restricts $\lambda$ to perfect squares, etc. 
\begin{table}[ht]
\begin{center}
\begin{tabular}{lll}
$p^1(4)=5$&$p^2(4)=2$&$p^{k\geq3}(4)=1$\\
&&\\
$4$&$2^2$&$1^k+1^k+1^k+1^k$\\
$3+1$&$1^2+1^2+1^2+1^2$&\\
$2+2$&&\\
$2+1+1$&&\\
$1+1+1+1$&&
\end{tabular}\
\end{center}
\end{table}\\
Independently, both Benfield, Paul, \& Roy \cite{benfield} and O'Sullivan \cite{[Sullivan]} proved that $p^k(n)$ is log-concave for all $k$ and sufficiently large $n$.
 \begin{theorem}\label{log-concavity}
     There exists an index $N_k$ such that for all $k$, the inequality $p^k(m)p^k(n)\geq p^k(m+n)$
     holds for $n\geq m \geq N_k$. 
 \end{theorem}
Conjecture \ref{p^k(n)_conj} gives instances $(m,n)$ where the inequality $p^2(m)p^2(n) \geq p^2(m+n)$ is equal and when it fails. The case $k=1$ was proven by Bessenrodt \& Ono \cite{[Bessenrodt]}.
\begin{conjecture}\label{p^k(n)_conj}
The square power partition function satisfies the inequality $p^2(m)p^2(n)\geq p^2(m+n)$ with instances of failure or equality as given in Table $2$:
\begin{table}[ht]
\begin{adjustbox}{width=1\textwidth}
\begin{tabular}{l|l}
Equality at $(m,n)$ with $n\geq m$&Failure at $(m,n)$ with $n\geq m$\\
\hline
$m=1$, $n=1, 2, 4, 5, 6, 9, 10, 13, 14, 18, 22$&$m=1$, $n\neq1, 2, 4, 5, 6, 9, 10, 13, 14, 18, 22$\\
$m=2, n=4,5,9,13$&$m=2$, $n\neq 1,4,5,9,13$\\
$m=3$, $n=4$&$m=3$, $n\neq1,2,4$\\
$m=4$, $n=5,6,7$&\\
$m=5$, $n=5,6,8,11,15$&$m=5$, $n=7$\\
$m=6$, $n=8,10,12,14$&$m=6$, $n=6,7,11,15$\\
$m=7$, $n=8,9,12,13,19$&$m=7$, $n=7,10,11,14,15$\\
\end{tabular}
\end{adjustbox}
\caption{Equality and Failure for $p^2(m)p^2(n)\geq p^2(m+n)$}
\end{table}
\end{conjecture}

What is the smallest index for which log-concavity fails? Define for $p^k(n)$, the smallest index $N_k$ such that for all $n>N_k$, the sequence ${p^k(n)}$ is log-concave. 
 \begin{conjecture}The smallest $N_k$ such that $p^k(n)$ is log-concave for all $n>N_k$ for $k~\leq~6$ is given by:
\begin{table}[ht]
\caption{Smallest $N_k$ for $k=2,\ldots,6$, OEIS: A346160 \cite{[A346160]}}
 \begin{center}
\begin{tabular}{c|cccccccc}
$p^k(n)$&$p^2(n)$&$p^3(n)$&$p^4(n)$&$p^5(n)$&$p^6(n)$\\
\hline
$N_k$&$1041$&$15655$&$637854$&$2507860$&$35577568$
\end{tabular}\
\end{center}
\end{table}
\end{conjecture}

\subsection{Fractional Partition Function} \ 
\newline
Let $\alpha$ be a real number and define the \textit{fractional partition function} $p_\alpha(n)$ by
    \[
    \sum_{n=0}^\infty p_\alpha(n)q^n = (q;q)_\infty^\alpha    \]
In 2019 Huat \& Chan \cite{[Huat]} studied the properties of $p_\alpha(n)$ where $\alpha$ is a rational number. The next year, Iskander, Jain, \& Talvola \cite{[Iskander]} proved that for any real $\alpha$ there exists an integer $N_{\alpha}$ such that the fractional partition function $p_\alpha(n)$ is log-concave $n>N_{\alpha}$. 
\begin{theorem}\label{fractional partition function}
Let $\alpha>0$. Then there is an index $M_\alpha\geq N_{\alpha}$ such that the sequence $p_\alpha(n)$ satisfies the inequality $p_\alpha(n)p_\alpha(m)>p_\alpha(n+m)$ for all $n,m>M_\alpha$ 
\end{theorem}

\subsection{The Multiset Partition Function}\ 
\newline
Let $k$ be a fixed positive integer and let $\mathcal{A}=\{a_i\}_{i=1}^\infty$ be a non-decreasing sequence of positive integers and define a restricted partition of $n$ as the sequence of positive integers
$\lambda_1,\lambda_2,\ldots,\lambda_j$ such that
$n = \lambda_1+\lambda_2+\ldots+\lambda_j$ where each $\lambda_i$ belongs to the multiset $\{a_1, a_2, . . . , a_k\}$. The \textit{multiset partition function} $p_\mathcal{A}(n, k)$ counts all
partitions of $n$ among the first $k$ elements of the multiset $\mathcal{A}$. The sequence $\{p_\mathcal{A}(n,k)\}$ was studied in 2022 by Gajdzica \cite{[Gajdzica]} who proved log-concavity under certain conditions of $n$ and $k$. This result was further improved by Gajdzica, Heim, \& Neuhauser \cite{gajdzica_heim_markus}.
\begin{theorem}\label{multiset partition function}
    There exists a sufficiently large index such that the multiset partition function satisfies the inequality $p_\mathcal{A}(n,k)p_\mathcal{A}(m,k) \geq p_\mathcal{A}(n+m,k)$ under the same conditions Gajdzica gives in his Theorem 1.1 \cite{[Gajdzica]}. 
\end{theorem}

\section{Log-polynomial Sequences}

In a more general setting, O'Sullivan \cite{[Sullivan2]} made interesting connections between the data of log-polynomial sequences and log-concavity/log-convexity. A sequence $a(n)$ is \textit{log-polynomial} of degree $m \geq 2$ with data $\{A(n),\kappa,\delta(n)\}$ for $\kappa=\pm1$ if there exists $\delta(n)>0$, $A(n)$,  and $g_k(n)$ for $k=3,4,\ldots,m$ such that
\[
\log{\left(\frac{a(n+j)}{a(n)}\right)}=A(n)j+\kappa\cdot\delta(n)^2j^2+\sum_{k=3}^mg_k(n)j^k+o\left(\delta(n)^{m+1}\right)
\]
for $j=1,2,\ldots,m+1$ where both $\delta(n)\rightarrow0$ and $g_k(n)=o\left(\delta(n)^k\right)$ as $n\rightarrow\infty$. O'Sullivan \cite{[Sullivan2]} showed that if $a(n)$ is log-polynomial for degrees $m\geq 2$, then $a(n)$ is log-concave for large $n$ whenever $\kappa=-1$ and is log-convex for large $n$ whenever $\kappa=1$. In particular, he showed \cite{[Sullivan2]} that if $a(n)=n^re^{t\cdot n^s}$ for real numbers $r$, $s$, and $t$ where $r$ and $t$ are not both zero, $0<s<2$, and $s(s-1)t\neq0$, then $a(n)$ is log-polynomial with data $\{A(n),\kappa,\delta(n)\}$ where 
\[
A(n)=-\frac{r}{n}+st\cdot n^{s-1}\quad \ \quad \kappa=\text{sgn}(s(s-1)t) \quad \ \quad \kappa\cdot\delta(n)^2=-\frac{r}{2n^2}+\frac{s(s-1)t}{2n^{2-s}}
\]
for all degrees $m \geq 2$.
\begin{theorem}
    If a sequence $a(n)$ is log-polynomial for degrees $m\geq 2$ with data $\{A(n),\kappa,\delta(n)\}$ and $0<s\leq 1$, then for sufficiently large integers $n$ and $m$, the sequence satisfies $a(n)a(m) \geq a(n+m)$ whenever $\kappa=-1$ and satisfies $a(n)a(m) \leq a(n+m)$ whenever $\kappa=1$ with finitely many instances of equality or failure.
\end{theorem}

\section{Sufficient But Not Necessary Condition for Partition Functions}
Although many restricted partition functions are log-concave (and hence abundant), there are restricted partition functions that are neither log-concave nor log-convex yet still appear to be abundant. Below is an empirical result for a particular restricted partition function.

\subsection{The \texorpdfstring{$m$}\ -ary Partition Function}\ 
\newline
The \textit{$m$-ary partition function} $b^m(n)$ counts the number of partitions of $n$ among powers of $m>1$. Asymptotics for the $m$-ary partitions functions were established by Kachi \cite{[Kachi]} in 2015. The most famous of the $m$-ary partition functions is the \textit{binary partition function} which counts the partitions of $n$ among powers of $2$. For example, $b^2(4)=4$, $b^3(4)=2$, $b^4(4)=2$ and for $m\geq5$, $b^m(4)=1$:
\begin{table}[ht]
\begin{center}
\begin{tabular}{cccc}
$b^2(4)=4$&$b^3(4)=2$&$b^4(4)=2$&$b^{m\geq5}(4)=1$\\
&&&\\
$2^2$&$3^1+3^0$&$4^1$&$m^0+m^0+m^0+m^0$\\
$2^1+2^1$&$3^0+3^0+3^0+3^0 $&$4^0+4^0+4^0+4^0$&\\
$2^1+2^0+2^0$&&&\\
$2^0+2^0+2^0+2^0$&&&
\end{tabular}\
\end{center}
\end{table}
\\
In 2019, Ganguly \cite{[Ganguly]} showed that $b^m(\alpha)$ is neither log-concave nor log-convex. Yet, for each $m$, the sequence $\{b^m(\alpha)\}_{\alpha\in\mathbb{N}}$ appears to be abundant with instances of equality or failure according to congruence classes for $\alpha<m$ and only finitely many instances of equality (and no instances of failure) for $\alpha\geq m$. A regular pattern appears for $m\geq6$.
\begin{conjecture}\label{m-ary partition}
For each $2\leq m\leq 5$, the sequence given by $\{b^m(\alpha)\}_{\alpha\in\mathbb{N}}$ satisfies the inequality $b^m(\alpha)b^m(\beta)~\geq~b^m(\alpha+\beta)$ with instances of equality or failure as given by Table \ref{mary_table}.
\begin{table}[ht]
\begin{adjustbox}{width=.75\textwidth}
\begin{tabular}{l|l|l}
$b^m(n)$&Equality at $(\alpha,\beta)$ with $\beta\geq \alpha$&Failure at $(\alpha,\beta)$ wtih $\beta\geq \alpha$\\
\hline
&$\alpha=1,\  \beta\equiv0\pmod{2}$&\\
$b^2(n)$&$\alpha=2$, $\beta=2, 3$&$\alpha=3 \quad \beta=3,5,7$\\
&$\alpha=3, \beta=9$&\\
\hline
&$\alpha=1,\ \beta\equiv0,1\pmod{3}$&$\alpha=1,\ \beta\equiv2\pmod{3}$\\
$b^3(n)$&$\alpha=2,\ \beta\equiv0\pmod{3}$&$\alpha=2,\ \beta\equiv1,2\pmod{3}$\\
&$\alpha=8,\ \beta=7,8$&$\alpha=4,\ \beta=8$\\
&&$\alpha=5,\ \beta=4,5,7,8$\\
\hline
&$\alpha=1,\ \beta\equiv0,1,2\pmod{4}$&$\alpha=1,\ \beta\equiv3\pmod{4}$\\
$b^4(n)$&$\alpha=2,\ \beta\equiv0,1\pmod{4}$&$\alpha=2,\ \beta\equiv2,3\pmod{4}$\\
&$\alpha=3,\ \beta\equiv0\pmod{4}$&$\alpha=3,\ \beta\equiv1,2,3\pmod{4}$\\
&$\alpha=5,\ \beta=7,11,15$&\\
&$\alpha =6,\ \beta=6,7,10,11,14,15$&\\
&$\alpha=7,\ \beta=7,9,10,11,13,14,15$&\\
\hline
&$\alpha=1,\ \beta\equiv0,1,2,3\pmod{5}$&$\alpha=1,\ \beta\equiv4\pmod{5}$\\
$b^5(n)$&$\alpha=2,\ \beta\equiv0,1,2\pmod{5}$&$\alpha=2,\ \beta\equiv3,4\pmod{5}$\\
&$\alpha=3,\ \beta\equiv0,1\pmod{5}$&$\alpha=3,\ \beta\equiv2,3,4\pmod{5}$\\
&$\alpha=4,\ \beta\equiv0\pmod{5}$&$\alpha=4,\ \beta\equiv1,2,3,4\pmod{5}$\\
&$7\leq\alpha\leq9$,\ $6\leq\beta\leq9$\\
\hline
\end{tabular}
\end{adjustbox}
\caption{Equality and Failure for $b^m(\alpha)b^m(\beta)\geq b^m(\alpha+\beta)$}
\label{mary_table}
\end{table}

\end{conjecture}
\section{Final Thoughts}
Log-concavity (or log-convexity) along with Condition \eqref{Condition 1} (or with Condition \eqref{Condition 2}) have been shown to be sufficient conditions for a sequence to be abundant (or deficient). Perhaps more interesting are sequences that are neither log-concave nor log-convex yet are still abundant or deficient. There are essentially seven types of sequences that are neither log-concave nor log-convex, according to whether the sequence adheres to Conditions \eqref{Condition 1} or \eqref{Condition 2} or neither, and whether it is abundant, deficient, or neither. 

\indent The first two types of sequences are those that are neither log-concave nor log-convex and satisfy neither Conditions \eqref{Condition 1} nor \eqref{Condition 2} yet are abundant or deficient. Example 1 presents a sequence that is neither log-concave nor log-convex, and satisfies neither Conditions \eqref{Condition 1} nor \eqref{Condition 2}, yet, the sequence is deficient for all indices $n>1$. Taking the reciprocal sequence of Example 3 reverses the inequalities in each property, producing a sequence that is neither log-concave nor log-convex, and satisfies neither Conditions \eqref{Condition 1} nor \eqref{Condition 2}, yet, the sequence is abundant for all indices $n>1$.

The next two types of sequences are those that are neither log-concave nor log-convex, but do satisfy Condition \eqref{Condition 1}, and are abundant or deficient. 
Taking the reciprocal the sequence found in Example 4 (the Fibonacci sequence) reverses the inequalities in each property, producing a sequence that is neither log-concave nor log-convex, but does satisfy Condition \eqref{Condition 1}, and is abundant.
\begin{question}
    Does there exist a sequence that is neither log-concave nor log-convex, yet satisfies Condition \eqref{Condition 1} and is deficient?
\end{question}

The next two types of sequences are those that are neither log-concave nor log-convex, but do satisfy Condition \eqref{Condition 2}, and are abundant or deficient. Example 4, the Fibonacci sequence, is neither log-concave nor log-convex, but does satisfy Condition \eqref{Condition 2}, and is deficient. 
\begin{question}
    Does there exist a sequence that is neither log-concave nor log-convex, yet satisfies Condition \eqref{Condition 2} and is abundant?
\end{question}

The last type of sequences are those that are neither log-concave nor log-convex, satisfy neither Conditions \eqref{Condition 1} nor \eqref{Condition 2}, and are neither abundant nor deficient. There are certainly such sequences; more investigation is required to categorize exactly which sequences satisfy these curious conditions. 

\section*{Data Sharing and Conflicts of Interest}
On behalf of all authors, the corresponding author states that there is no conflict of interest. Data sharing not applicable to this article as no datasets were generated or analysed during the current study. 

\section*{Acknowledgements}
The author would like to extend sincere gratitude to the anonymous referee, whose keen observations and criticism throughout the review process has greatly improved the quality of this paper.


\begin{thebibliography}{10}
	
	\bibitem{[Ahmiaab_Belbachirb_Belkhirb]}
	Moussa Ahmia, Hac{\`e}ne Belbachir, and Amine Belkhir.
	\newblock The log-concavity and log-convexity properties associated to
	{H}yper{P}ell and {H}yper{P}ell-{L}ucas sequences.
	\newblock {\em Annales Mathematicae et Informaticae}, 43:3--12, 2014.
	
	\bibitem{[Alanzi]}
	Abdulaziz~A. Alanazi, Stephen~M. Gagola, and Augustine~O. Munagi.
	\newblock Combinatorial proof of a partition inequality of {B}essenrodt-{O}no.
	\newblock {\em Annals of Combinatorics}, 21:331--337, 2017.
	
	\bibitem{[Asai]}
	Nobuhiro Asai, Izumi Kubo, and Hui-Hsiung Kuo.
	\newblock Bell numbers, log-concavity, and log-convexity.
	\newblock {\em Acta Applicandae Mathematica}, 63(1-3):79--87, 2000.
	
	\bibitem{[Beckwith]}
	Olivia Beckwith and Christine Bessenrodt.
	\newblock Multiplicative properties of the number of k-regular partitions.
	\newblock {\em Annals of Combinatorics}, 20:231--250, 2016.
	
	\bibitem{benfield}
	Brennan Benfield, Madhumita Paul, and Arindam Roy.
	\newblock Tur{\'a}n inequalities for k-th power partition functions.
	\newblock {\em Journal of Mathematical Analysis and Applications},
	529(1):127678, 2024.
	
	\bibitem{BessOno}
	Christine Bessenrodt and Ken Ono.
	\newblock Maximal multiplicative properties of partitions.
	\newblock {\em Ann. Comb.}, 20(1):59--64, 2016.
	
	\bibitem{BKRT}
	Kathrin Bringmann, Ben Kane, Larry Rolen, and Zack Tripp.
	\newblock Fractional partitions and conjectures of {C}hern-{F}u-{T}ang and
	{H}eim-{N}euhauser.
	\newblock {\em Trans. Amer. Math. Soc. Ser. B}, 8:615--634, 2021.
	
	\bibitem{[Huat]}
	Heng~Huat Chan and Liuquan Wang.
	\newblock Fractional powers of the generating function for the partition
	function.
	\newblock {\em Acta Arithmetica}, 187:59--80, 2018.
	
	\bibitem{[Chen]}
	W.~Y.~C. Chen.
	\newblock Recent developments on the log-concavity and q-log-concavity of
	combinatorial polynomials.
	\newblock {\em FPSAC 2010 Conference Talk Slides}, 2010.
	
	\bibitem{[Chern]}
	Shane Chern, Shishuo Fu, and Dazhao Tang.
	\newblock Some inequalities for k-colored partition functions.
	\newblock {\em The Ramanujan Journal}, 46:713--725, 2018.
	
	\bibitem{[Craig]}
	William Craig and Anna Pun.
	\newblock A note on the higher order {T}ur{\'a}n inequalities for k-regular
	partitions.
	\newblock {\em Research in Number Theory}, 7(1):5, 2021.
	
	\bibitem{[Salvo]}
	Stephen DeSalvo and Igor Pak.
	\newblock Log-concavity of the partition function.
	\newblock {\em The Ramanujan Journal}, 38:61--73, 2015.
	
	\bibitem{[Dong]}
	Janet J.~W. Dong and Kathy~Q. Ji.
	\newblock Higher order {T}ur\'an inequalities for the distinct partition
	function.
	\newblock {\em J. Number Theory}, 260:71--102, 2024.
	
	\bibitem{[Engel]}
	Benjamin Engel.
	\newblock Log-concavity of the overpartition function.
	\newblock {\em The Ramanujan Journal}, 43:229--241, 2017.
	
	\bibitem{[Euler]}
	Leonhard Euler.
	\newblock {\em Introductio in analysin infinitorum}, volume~2.
	\newblock Apud Marcum-Michaelem Bousquet \& Socios, 1748.
	
	\bibitem{[Gajdzica]}
	Krystian Gajdzica.
	\newblock Log-concavity of the restricted partition function $p_\mathcal{A}(n,
	k)$ and the new {B}essenrodt-{O}no type inequality.
	\newblock {\em Journal of Number Theory}, 251:31--65, 2023.
	
	\bibitem{gajdzica_heim_markus}
	Krystian Gajdzica, Bernhard Heim, and Markus Neuhauser.
	\newblock Polynomization of the {B}essenrodt-{O}no type inequalities for
	{A}-partition functions.
	\newblock {\em Annals of Combinatorics}, pages 1--23, 2024.
	
	\bibitem{Gajdzica_Miska_Ulas}
	Krystian Gajdzica, Piotr Miska, and Maciej Ulas.
	\newblock On general approach to {B}essenrodt-{O}no type inequalities and
	log-concavity property.
	\newblock {\em Annals of Combinatorics}, 2024.
	
	\bibitem{[Ganguly]}
	Soumendra Ganguly.
	\newblock {\em $m$-ary Partitions}.
	\newblock PhD thesis, TigerPrints, 2019.
	
	\bibitem{[Gedefa_2]}
	Fekadu~Tolessa Gedefa.
	\newblock Log-concavity of centered polygonal figurate number sequences.
	\newblock {\em Open Access Library Journal}, 3(06):1, 2016.
	
	\bibitem{[Gedefa_1]}
	Fekadu~Tolessa Gedefa.
	\newblock On the log-concavity of polygonal figurate number sequences.
	\newblock {\em arXiv preprint arXiv:2006.05286}, 2020.
	
	\bibitem{heimmarkus}
	Bernhard Heim and Markus Neuhauser.
	\newblock Sequence with inequalities.
	\newblock {\em https://arxiv.org/abs/2408.00319}.
	
	\bibitem{[Heim]}
	Bernhard Heim, Markus Neuhauser, and Robert Tr{\"o}ger.
	\newblock Inequalities for plane partitions.
	\newblock {\em Annals of Combinatorics}, 27(1):87--108, 2023.
	
	\bibitem{[Honsberger]}
	Ross Honsberger.
	\newblock A second look at the {F}ibonacci and {L}ucas numbers.
	\newblock {\em Mathematical gems III. Washington, DC: Math. Assoc. Amer.},
	1985.
	
	\bibitem{[Iskander]}
	Jonas Iskander, Vanshika Jain, and Victoria Talvola.
	\newblock Exact formulae for the fractional partition functions.
	\newblock {\em Research in Number Theory}, 6(2):20, 2020.
	
	\bibitem{[Kachi]}
	Yasuyuki Kachi and Pavlos Tzermias.
	\newblock On the $m$-ary partition numbers.
	\newblock {\em Algebra and Discrete Mathematics}, 2015.
	
	\bibitem{[Kim]}
	Dongsu Kim and Ae~Ja Yee.
	\newblock A note on partitions into distinct parts and odd parts.
	\newblock {\em The Ramanujan Journal}, 3(2):227--231, 1999.
	
	\bibitem{[Lehmer]}
	Derrick~H. Lehmer.
	\newblock On the partition of numbers into squares.
	\newblock {\em The American Mathematical Monthly}, 55(8):476--481, 1948.
	
	\bibitem{liu}
	Edward Y.~S. Liu and Helen W.~J. Zhang.
	\newblock Inequalities for the overpartition function.
	\newblock {\em The {R}amanujan Journal}, 54:485--509, 2021.
	
	\bibitem{[Dawsey]}
	Madeline Locus~Dawsey and Riad Masri.
	\newblock Effective bounds for the {A}ndrews spt-function.
	\newblock {\em Forum Mathematicum}, 31(3):743--767, 2019.
	
	\bibitem{[Merca]}
	Mircea Merca.
	\newblock On the partitions into distinct parts and odd parts.
	\newblock {\em Quaestiones Mathematicae}, 44(8):1095--1105, 2021.
	
	\bibitem{[Nicolas]}
	Jean-Louis Nicolas.
	\newblock Sur les entiers $ n $ pour lesquels il y a beaucoup de groupes
	{A}b{\'e}liens d’ordre $ n$.
	\newblock {\em Annales de l'institut Fourier}, 28(4):1--16, 1978.
	
	\bibitem{[Bessenrodt]}
	Ken Ono.
	\newblock Distribution of the partition function modulo $m$.
	\newblock {\em Annals of Mathematics}, 151(1):293--307, 2000.
	
	\bibitem{[Onoetal]}
	Ken Ono, Sudhir Pujahari, and Larry Rolen.
	\newblock {T}ur{\'a}n inequalities for the plane partition function.
	\newblock {\em Advances in Mathematics}, 409:108692, 2022.
	
	\bibitem{[Sullivan]}
	Cormac O'Sullivan.
	\newblock Detailed asymptotic expansions for partitions into powers.
	\newblock {\em Int. J. Number Theory}, 19(9):2163--2196, 2023.
	
	\bibitem{[Sullivan2]}
	Cormac O’Sullivan.
	\newblock Limits of {J}ensen polynomials for partitions and other sequences.
	\newblock {\em Monatshefte f{\"u}r Mathematik}, 199(1):203--230, 2022.
	
	\bibitem{[SIBUY]}
	Masaaki Sibuya.
	\newblock Log-concavity of {S}tirling numbers and unimodality of {S}tirling
	distributions.
	\newblock {\em Annals of the Institute of Statistical Mathematics},
	40:693--714, 1988.
	
	\bibitem{[000009]}
	Neil J.~A. Sloane and The OEIS~Foundation Inc.
	\newblock The on-line encyclopedia of integer sequences {A}000009.
	\newblock https://oeis.org/A000009, 2023.
	
	\bibitem{[A015128]}
	Neil J.~A. Sloane and The OEIS~Foundation Inc.
	\newblock The on-line encyclopedia of integer sequences {A}015128.
	\newblock https://oeis.org/A015128, 2023.
	
	\bibitem{[A346160]}
	Neil J.~A. Sloane and The OEIS~Foundation Inc.
	\newblock The on-line encyclopedia of integer sequences {A}346160.
	\newblock https://oeis.org/A346160, 2023.
	
	\bibitem{[Spicey]}
	Michael~Z. Spivey.
	\newblock {F}ibonacci identities via the determinant sum property.
	\newblock {\em The College Mathematics Journal}, 37(4):286--289, 2006.
	
	\bibitem{[Zheng]}
	Li-Na Zheng, Rui Liu, and Feng-Zhen Zhao.
	\newblock On the log-concavity of the {H}yper{F}ibonacci numbers and the
	{H}yper{L}ucas numbers.
	\newblock {\em Journal of Integer Sequences}, 17(1):14--1, 2014.
	
\end{thebibliography}

\end{document}